\documentclass{amsart}
\usepackage{graphicx}
\numberwithin{equation}{section}

\begin{document}
\title[A New Optimum Tuning Method in Time-Delay Systems]{A New Optimum Tuning Method of PI Controllers in First-Order Time-Delay Systems}


\author{Gianpasquale Martelli}
\address{Via Domenico da Vespolate 8 28079 Vespolate Italy}
\curraddr{}
\email{gianpasqualemartelli@libero.it}
\thanks{}



\date{May 16, 2007}

\dedicatory{}

\begin{abstract}
In this paper a new optimum tuning method of PI controllers in first-order time-delay systems, based on the deadbeat response to a step setpoint variation, is presented. The deadbeat performance, already studied for the plants without delay, consists of a fast achievement of the steady state with an overshoot included in a narrow band and with minimum rise and settling times. In the proposed method the rise and settling times are both replaced with the integral of the squared error and the constraint, which requires a controller output less than a preset value, is added. The transient behaviour is strictly evaluated by means of the  analytical solutions of the relative differential difference equations, obtained by the method of steps for the first time. Moreover the proposed tuning method is compared with other three methods, selected among the most used, which are the Ziegler-Nichols time-domain and frequency-domain  and the Zhuang and Atherton ISTE ones. Finally tuning charts, having as coordinates the two PI controller parameters and provided with the borderlines of the stability region, are introduced.
\end{abstract}

\maketitle

\section{Introduction}
The implementation of an optimum setting of the proportional integral derivative (PID) controllers , especially in time-delay systems, is one of the most attractive goal of the control designers and it has been investigated from more than 60 years. It seems unbelievable  that today an original and interesting word can be said, but this happens in the present paper for the PI controllers thanks to a new mathematical tool, which is the explicit solution of differential difference equations detailed in \cite{bib1}. 

Before the enunciation of the proposed method, let us mention some of the most significative existing ones. The first attempt has been performed by Ziegler and Nichols, who presented two design methods, one time-domain and one frequency-domain \cite{bib2}. In both methods the amplitude of the decay ratio, which is the ratio between two consecutive maxima of the control error during a load disturbance, has been fixed equal to one quarter.
The Ziegler-Nichols tuning rules have the advantage of being easy to use, but give closed-loop systems with very poor damping. Modified expressions of the controller parameters, which improve the damping, have been proposed by Chien, Hrones and Reswick \cite{bib3} and by Cohen and Coon \cite{bib4}.
A second method, named Kappa-Tau \cite{bib5}, assumes the maximum sensivity as design parameter and applies the dominant pole design to typical process control models.
A third method is based on the minimization of an integral performance criterion, whose the most well known is the integral of  the squared error ($ISE$).  Zhuan and Atherton \cite{bib6} have approximated the exponential term of the first-order time-delay plant model using a 3/3 order Pad\'e and have minimized the $ISE$, $ISTE$ and $IST^{2}E$ performance criteria for setpoint and load changes.
An other method, which deserves to be cited, is the deadbeat response one, studied for processes without time delay and defined in para. 10.2.5 of \cite{bib7} pag. 171 as follows:
\begin{itemize}
\item Zero steady-state error
\item Minimum rise time and settling time
\item $0.1 \% <=$ controlled variable overshoot $<= 2 \%$
\item Controlled variable undershoot $<= 2 \%$
\end{itemize}

The proposed method, which is a modification of the deadbeat response one, can be defined concisely, but exhaustively, as follows: it minimizes the integral squared error ($ISE$), during a closed-loop step setpoint change, with the controlled variable overshoot ($PO_{y}$) equal to a preset value (for example $ 1.05 \%$) and the controller output overshoot ($PO_{v}$) lower than a preset value (for example $ 10 \%$).
Only two performances indexes, $ISE$ and $PO_{y}$ or $PO_{v}$, are involved in this tuning optimization and therefore a PI controller, provided with two parameters, is adequate i.e. a PID controller is not required.

A significative advantage of this method is the avoidance of the controller windup without any additional sophisticated scheme.
It is known that one of the nonlinear issues of any controller is the actuator saturation, since from the practical point of view the actuator transfer function is bonded and the controller output may be out of the rated band. When the actuator saturates, the feedback loop is broken and the actuator output is constant and then independent of the error signal. In this case a large overshoot of the controlled variable may occur, a long time may be needed for the recovery of the system and the actuator may bounce several times between the two bond values.
If the rated value of the controller output is set lower than the extreme point of the linear part of the actuator transfer function and if the controller is tuned in such a way that the known maximum overshoot of the controller output never reaches this extreme point, none windup occurs.

The feedback system taken into consideration is represented by the block diagram of Fig. 1 and the transfer functions are mathematically described by
\begin{equation}\label{eq:1.1}
G_{p}(s)=K  \frac{e^{-L \, s}}{1+T_{p}s}\
\end{equation}
\begin{equation}\label{eq:1.2}
G_{c}(s)=K_{p}+ \frac{K_{i}}{s}\ 
\end{equation}
or
\begin{equation}\label{eq:1.3}
G_{c}(s)=K_{p} \left(1 +\frac{1}{T_{i}s}\  \right)
\end{equation}
where $K$ is the plant steady-state gain, $T_{p}$ the positive plant time constant, $L$ the plant time delay and $K_{p}$, $K_{i}$ and $T_{i}$  the parameters of the PI controller.
\begin{figure}[htbp]
\centering
\includegraphics[bb= -25 -51 181 16]{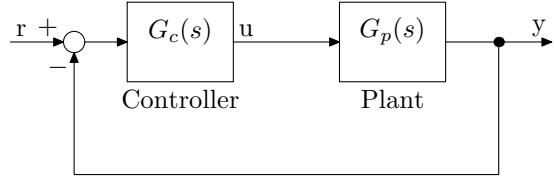}
\caption{Feedback control system}
\end{figure}

In order to get equations independent of the real values of the parameters, the normalized time $t$ and dimensionless parameters are introduced as follows:
\begin{itemize}
\item Systems without delay \\
$t$ referred to the plant time costant $T_{p}$ and  $h=K \,K_{p}$, $t_{i} = T_{i}/T_{p}$.
\item Systems with delay \\
$t$ referred to the plant time delay $L$ and  $t_{p}=T_{p}/L$, $h=K \,K_{p}$, $h_{i} = K \,K_{i}L$.
\end{itemize}

The open loop transfer functions for this system without time delay $F_{1}( \sigma)$  ($\sigma=T_{p}s $ in (\ref{eq:1.1}) and (\ref{eq:1.3})) and with time delay  $F_{2}(\sigma)$ ($\sigma=L \,s$ in (\ref{eq:1.1}) and (\ref{eq:1.2})) are given respectively by
\begin{equation}\label{eq:1.4}
F_{1}(\sigma)= \frac{h(1+t_{i} \sigma)}{t_{i} \sigma(1+\sigma)}\ 
\end{equation}
\begin{equation}\label{eq:1.5}
F_{2}(\sigma)=  \frac{h_{i}+h \, \sigma}{\sigma \left(1+t_{p} \sigma \right)}\ e^{-\sigma}
\end{equation}

The tuning points will be positioned in tuning charts, whose coordinates are $h$ and $t_{i}$ for systems without delay and $h$ and $h_{i}$ for systems with delay. The stability region borderlines and the curves of constant phase margin will be also plotted in these charts and therefore the stability of each selected tuning point may be easily verified.

In the next Sections the proposed tuning method will be detailed, applied to both systems, without and with time delay, and finally compared with the two Ziegler-Nichols and with the Zhuang-Atherton $ISTE$ for setpoint change ones.

\section{Systems without time delay}
For $h+1>0$ the system is stable and hence all the tuning points of the quadrant $h>0$ and $t_{i}>0$ are stable. The equations of the curves  $\Gamma_{y}$, $\Gamma_{v}$ and $\Gamma_{i}$, relative respectively to the controlled variable overshoot $PO_{y}$,  the controller output overshoot $PO_{v}$  and  the integral of the square error $ISE$, are given by  (\ref{eq:A.8}), (\ref{eq:A.9}) and (\ref{eq:A.10}). Equating (\ref{eq:A.4}) to zero, one obtains for the damping borderline
\begin{displaymath}
t_{i} = \frac{4 \,h}{(1+h)^{2}}\
\end{displaymath}
The response is overdamped  if the selected point is above this line and underdamped if below.
In Fig. 2 the following curves are plotted:
\begin{itemize}
\item $\Gamma_{d}$: damping borderline
\item $\Gamma_{v}$:  $PO_{v}=0.1$
\item $\Gamma_{y1}$, $\Gamma_{y2}$, $\Gamma_{y3}$: respectively $PO_{y}=0.1 \%$, $PO_{y}=1.05 \%$, $PO_{y}=2 \%$
\item $\Gamma_{i1}$, $\Gamma_{i2}$, $\Gamma_{i3}$: respectively $ISE=1.0$, $ISE=1.2$, $ISE=1.4$. 
\end{itemize}
The optimum tuning point, named $B$ and corresponding to the minimum value of the $ISE$ for $PO_{y}=1.05 \%$ and $PO_{v}<=0.1$, is the intersection between $\Gamma_{y2}$ and $\Gamma_{v}$, since the $ISE$ always decreases when $PO_{v}$ increases for a given value of $PO_{y}$.
\begin{figure}[htbp]
\centering
\includegraphics[bb= -28 -55 318 254]{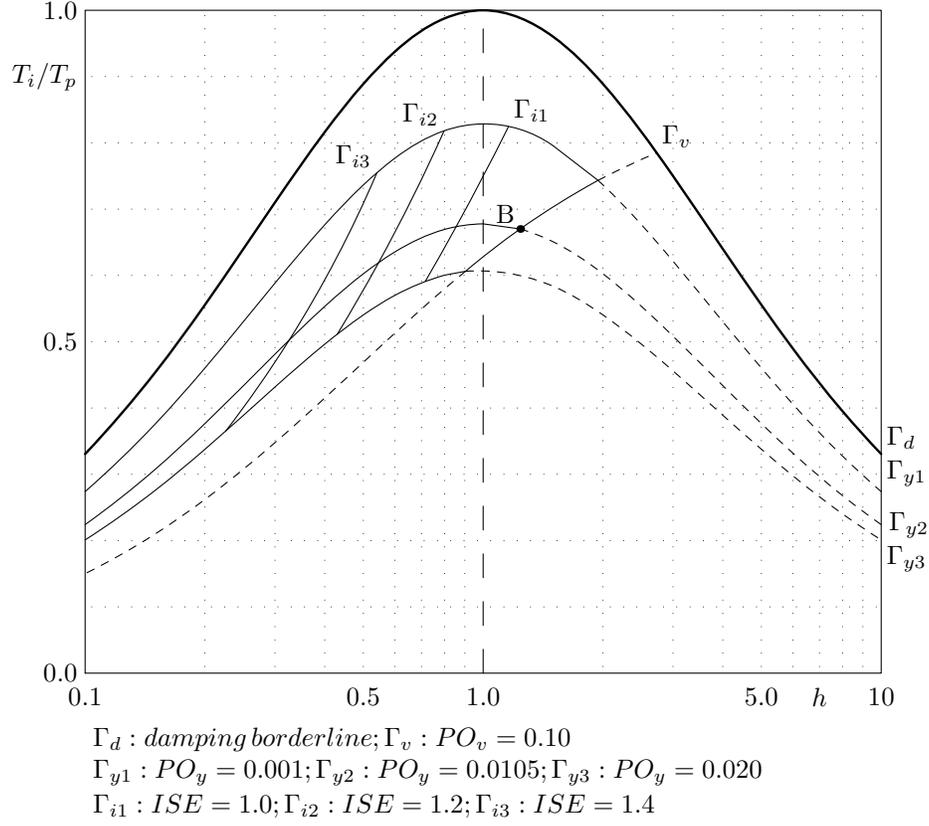}
\caption{Tuning chart of systems without time delay}
\end{figure}

\section{Systems with time delay}
\subsection{Stability region}
The stability zones are defined in \cite{bib8} and \cite{bib9} as follows:
\begin{enumerate}
\item The parameter $h$ must be included in the range from zero to $h_{u}$ given by
\begin{equation}\label{eq:3.1}
h_{u} = - \cos(z_{a})+t_{p}z_{a}sin(z_{a})
\end{equation}
where $z_{a}$ is the first positive solution of
\begin{displaymath}
tan(z_{a}) = - \frac{t_{p}}{1+t_{p}}\ z_{a}
\end{displaymath}
\item The parameter $h_{i}$, for a given value of $h$, must satisfy the following inequalities
\begin{equation}\label{eq:3.3}
\delta_{r}(z_{1})<0 \, and \, \delta_{r}(z_{2})>0 \, and \, h_{i}>0
\end{equation}
 where
 \begin{displaymath}
\delta_{r}(z)=h_{i} - {z} sin(z) - t_{p}z^{2} cos(z)
\end{displaymath}
and $z_{1}$ and $z_{2}$ are the first two positive roots of
\begin{displaymath}
h+cos(z)-t_{p}{z} sin(z)=0
\end{displaymath}
\end{enumerate}
The stability region borderline $\Gamma_{s}$ is plotted in Fig. 3 for $t_{p}=0.55$ and in Fig. 4 for $t_{p}=2.5$. 

It is also convenient to consider  the stability phase margin (PM), whose expression, deduced  from (\ref{eq:1.5}), is
\begin{displaymath}
\tan(z_{b}+PM)=-\frac{h_{i}+z_{b}^{2}h \,t_{p}}{z_{b}(h-h_{i}t_{p})}\
\end{displaymath}
where $z_{b}$ is the solution of
\begin{displaymath}
h^{2}+ \frac{h_{i}^{2}}{z_{b}^{2}}\ = 1+ t_{p}^{2}z_{b}^{2}
\end{displaymath}
The curves relative to $PM=30^{\circ}$, $PM=45^{\circ}$ and $PM=60 ^{\circ}$ are also included in Fig. 3 and Fig. 4. 

\subsection{Tuning}
The performance indexes $PO_{y}$, $PO_{v}$ and $ISE$ are evaluated in a period equal to seven times the process time delay, which is considered enough to obtain sound results, in order to reduce the computing complexity to an acceptable level. Therefore, for given values of $h$ and $h_{i}$, the controlled variable $y(t_{n})$ and afterwards the controller output $u(t_{n})$ multiplied by $K$, named $v(t_{n})$, are  calculated in 701 equally spaced points from $t=0$ to $t=7$ according to (\ref{eq:B.3}) and (\ref{eq:B.4}). The  opposite of the minimum negative value of $y(t_{n})$ is assumed  as $PO_{y}$ and the analogous value of $v(t_{n})$ as $PO_{v}$. The $ISE$ is computed with the triangle rule by means of  (\ref{eq:B.5}).  
In both Fig. 3, valid for $t_{p}=0.55$, and  Fig. 4, valid for $t_{p}=2.5$, the  curves $\Gamma_{y}$ for  $PO_{y}=0.0105$ and $\Gamma_{v}$ for  $PO_{v}=0.1$ are added.
The optimum tuning point, named $B_{4}$, has been determined as follows:
\begin{itemize}
\item let $C_{y}$ and $C_{v}$ denote respectively each point lying on the curves $\Gamma_{y}$ and $\Gamma_{v}$, corresponding to the 50 equally spaced values of $h$ from $h=0$ to the maximum allowed by the stability
\item compute the $ISE$ in  $C_{y}$, if $\Gamma_{y}$ is below $\Gamma_{v}$, or in $C_{v}$, if above, since both $PO_{y}$ and $PO_{v}$ increase when, for a given $h$, $h_{i}$ increases
\item $B_{4}$ is the point having the minimum value of the $ISE$
\end{itemize}
It has been found that for  $t_{p}=0.55$ (Fig. 3) $B_{4}$ lies on the curve $\Gamma_{y}$ and for $t_{p}=2.5$ (Fig. 4) on $\Gamma_{v}$.

\begin{figure}[htbp]
\centering
\includegraphics[bb= -20 -65 303 244]{martelliLf3.eps}
\caption{Tuning chart of systems with time delay- $t_{p}=0.55$}
\end{figure}

\begin{figure}[htbp]
\centering
\includegraphics[bb= -19 -65 303 204]{martelliLf4.eps}
\caption{Tuning chart of systems with time delay - $t_{p}=2.5$}
\end{figure}

\subsection{Comparison}
The proposed optimum tuning is compared with the Ziegler-Nichols time-domain and frequency domain \cite{bib2} and with Zhuang-Atherton $ISTE$ for setpoint change one \cite{bib6}.  Taking into consideration the equalities $a =K/t_{p}$,  $K_{u}K=(1+(z_{u}t_{p})^{2})^{0.5}$, $T_{u}/L=2 \pi /z_{u}$, the parameters are given by:
\begin{enumerate}
\item Ziegler-Nichols time-domain
\begin{displaymath}
K_{p}= \frac{0.9}{a}\ = \frac{0.9}{K}\ t_{p}; \quad T_{i}/L=3 
\end{displaymath}  
\item Ziegler-Nichols frequency-domain
\begin{displaymath}
K_{p}= 0.4 \, K_{u}= \frac{0.4}{K} \left(1+z_{u}^{2}t_{p}^{2} \right)^{0.5}; \quad  T_{i}/L=0.8 \, T_{u}/L= 0.8 \,2 \pi/z_{u}
\end{displaymath}  
where $z_{u}$ is the first positive root of $tan(z_{u}) = -z_{u}t_{p}$ 
\item Zhuang-Atherton $ISTE$ setpoint
\begin{enumerate}
\item $0.5<=tp<=0.9$
\begin{displaymath}
K_{p}= \frac{0.786}{K}\ (1/t_{p})^{-0.559}; \quad T_{i}/L= t_{p}/(0.883-0.158/t_{p}) 
\end{displaymath} 
\item $1<=tp<=10$
\begin{displaymath}
K_{p}=\frac{0.712}{K}\ (1/t_{p})^{-0.921}; \quad T_{i}/L=t_{p}/(0.968-0.247/t_{p})
\end{displaymath}   
\end{enumerate}
\end{enumerate}
Since $h=K \, K_{p}$ and $h_{i}=h \, L/T_{i}$, it is possible to calculate the optimum tuning points $B_{1}$, $B_{2}$ and $B_{3}$, also plotted in Fig. 3 and Fig. 4, for these three tuning rules.

The optimum controllers parameters $h$ and $h_{i}$ are listed in Table 1 whereas the performance indexes $PO_{y}$, $PO_{v}$ and $ISE$ are in Table 2 for some values of $t_{p}$, from $t_{p}=0.1$ to $t_{p}=10$. From Table 2 it follows that:
\begin{enumerate}
\item  the controlled variable response is overdamped if $t_{p}>=0.85$ for the proposed method and always for the others ones. 
\item  the controller output response is overdamped, if $t_{p}<=1$,  for the two Ziegler-Nichols and Zhuang-Atherton methods. Moreover the values of $PO_{v}$ for $t_{p}>1$ are unacceptable for the first two methods since higher than $0.10$, but good for the third one since lower than $0.06$. 
\item  neglecting the Ziegler-Nichols methods for $t_{p}>1$ in accordance with the previous consideration, the proposed method allows the best values of the $ISE$.
\end{enumerate}

\begin{table}[htbp]
\caption{Optimum controllers parameters}
\label{Tab1}
\centering
\begin{tabular}{|r @{.} l|r @{.} l|r @{.} l|r @{.} l|r @{.} l|r @{.} l|r @{.} l|r @{.} l|r @{.} l|}
\hline
\multicolumn{2}{|c|}{$t_{p}$} &
\multicolumn{8}{|c|}{$h$} &
\multicolumn{8}{|c|}{$h_{i}$} \\
\cline{3-18}
\multicolumn{2}{|c|}{} &
\multicolumn{2}{|c|}{Z-N} &
\multicolumn{2}{|c|}{Z-N} &
\multicolumn{2}{|c|}{Z-A} &
\multicolumn{2}{|c|}{Prop.} &
\multicolumn{2}{|c|}{Z-N} &
\multicolumn{2}{|c|}{Z-N} &
\multicolumn{2}{|c|}{Z-A} &
\multicolumn{2}{|c|}{Prop.} \\
\multicolumn{2}{|c|}{} &
\multicolumn{2}{|c|}{time} &
\multicolumn{2}{|c|}{freq.} &
\multicolumn{2}{|c|}{ISTE} &
\multicolumn{2}{|c|}{} &
\multicolumn{2}{|c|}{time} &
\multicolumn{2}{|c|}{freq.} &
\multicolumn{2}{|c|}{ISTE} &
\multicolumn{2}{|c|}{} \\
\hline
 0&10 & 0&090 & 0&416 & \multicolumn{2}{|c|}{} & 0&45 
      & 0&030 & 0&237 & \multicolumn{2}{|c|}{} & 0&787  \\
\hline
 0&25 & 0&225 & 0&475 & \multicolumn{2}{|c|}{} & 0&50 
      & 0&075 & 0&243 & \multicolumn{2}{|c|}{} & 0&738  \\
\hline 
 0&40 & 0&360 & 0&552 & \multicolumn{2}{|c|}{} & 0&60 
      & 0&120 & 0&262 & \multicolumn{2}{|c|}{} & 0&724  \\
\hline
 0&55 & 0&495 & 0&636 & 0&563 & 0&70 
      & 0&165 & 0&285 & 0&609 & 0&737  \\
\hline
 0&70 & 0&630 & 0&724 & 0&644 & 0&92 
      & 0&210 & 0&311 & 0&605 & 0&763 \\
\hline
 0&85 & 0&765 & 0&814 & 0&718 & 1&10  
      & 0&255 & 0&338 & 0&589 & 0&766  \\
\hline
 1&00 & 0&900 & 0&905 & 0&786 & 1&15  
      & 0&300 & 0&365 & 0&570 & 0&744  \\
\hline
 2&50 & 2&250 & 1&835 & 1&656 & 2&10
      & 0&750 & 0&654 & 0&576 & 0&682  \\
\hline
 4&00 & 3&600 & 2&774 & 2&553 & 3&00 
      & 1&200 & 0&947 & 0&578 & 0&654  \\
\hline
 5&50 & 4&950 & 3&715 & 3&423 & 3&80  
      & 1&650 & 1&241 & 0&574 & 0&633  \\
\hline
 7&00 & 6&300 & 4&656 & 4&274 & 4&75 
      & 2&100 & 1&535 & 0&569 & 0&628  \\
\hline
 8&50 & 7&650 & 5&598 & 5&111 & 6&00  
      & 2&550 & 1&829 & 0&564 & 0&640  \\
\hline
 10&00 & 9&000 & 6&540 & 5&936 & 6&65 
       & 3&000 & 2&123 & 0&560 & 0&622 \\
\hline
\end{tabular}
\end{table}

\begin{table}[htbp]
\caption{Optimum controllers performance indexes}
\label{Tab2}
\centering
\begin{tabular}{|r @{.} l|r @{.} l|r @{.} l|r @{.} l|r @{.} l|r @{.} l|r @{.} l|r @{.} l|r @{.} l|r @{.} l|}
\hline
\multicolumn{2}{|c|}{$t_{p}$} &
\multicolumn{2}{|c|}{$PO_{y}$} &
\multicolumn{8}{|c|}{$PO_{v}$} &
\multicolumn{8}{|c|}{$ISE$} \\
\cline{3-20}
\multicolumn{2}{|c|}{} &
\multicolumn{2}{|c|}{Prop.} &
\multicolumn{2}{|c|}{Z-N} &
\multicolumn{2}{|c|}{Z-N} &
\multicolumn{2}{|c|}{Z-A} &
\multicolumn{2}{|c|}{Prop.} &
\multicolumn{2}{|c|}{Z-N} &
\multicolumn{2}{|c|}{Z-N} &
\multicolumn{2}{|c|}{Z-A} &
\multicolumn{2}{|c|}{Prop.} \\
\multicolumn{2}{|c|}{} &
\multicolumn{2}{|c|}{} &
\multicolumn{2}{|c|}{time} &
\multicolumn{2}{|c|}{freq.} &
\multicolumn{2}{|c|}{ISTE} &
\multicolumn{2}{|c|}{} &
\multicolumn{2}{|c|}{time} &
\multicolumn{2}{|c|}{freq.} &
\multicolumn{2}{|c|}{ISTE} &
\multicolumn{2}{|c|}{} \\
\hline
 0&10 & 0&010 & \multicolumn{2}{|c|}{} & \multicolumn{2}{|c|}{}&
        \multicolumn{2}{|c|}{} & 0&014 & 6&095 & 3&186 & 
        \multicolumn{2}{|c|}{} & 1&524  \\
\hline
 0&25 & 0&010 & \multicolumn{2}{|c|}{} & \multicolumn{2}{|c|}{}&
        \multicolumn{2}{|c|}{} & 0&029 & 5&192 & 3&266 & 
        \multicolumn{2}{|c|}{} & 1&674  \\
\hline 
 0&40 & 0&010 & \multicolumn{2}{|c|}{} & \multicolumn{2}{|c|}{}&
        \multicolumn{2}{|c|}{} & 0&044 & 4&605 & 3&261 & 
        \multicolumn{2}{|c|}{} & 1&788  \\
\hline
 0&55 & 0&010 & \multicolumn{2}{|c|}{} & \multicolumn{2}{|c|}{}&
        \multicolumn{2}{|c|}{} & 0&086 & 4&193 & 3&229 & 
        1&998 & 1&869  \\
\hline
 0&70 & 0&010 & \multicolumn{2}{|c|}{} & \multicolumn{2}{|c|}{}&
        \multicolumn{2}{|c|}{} & 0&099 & 3&896 & 3&191 & 
        2&098 & 1&945  \\
\hline
 0&85 & \multicolumn{2}{|c|}{} & \multicolumn{2}{|c|}{} &  
        \multicolumn{2}{|c|}{} & \multicolumn{2}{|c|}{} &
         0&100 & 3&674 & 3&154 & 2&212 & 2&037  \\
\hline
 1&00 & \multicolumn{2}{|c|}{} & \multicolumn{2}{|c|}{} & 
        \multicolumn{2}{|c|}{} & \multicolumn{2}{|c|}{} &
         0&100 & 3&502 & 3&119 & 2&333 & 2&129  \\
\hline
 2&50 & \multicolumn{2}{|c|}{} & 0&177 & 0&112 & 0&032 & 0&100 &
        2&822 & 2&925 & 3&089 & 2&939  \\
\hline
 4&00 & \multicolumn{2}{|c|}{} & 0&849 & 0&609 & 0&049 & 0&100 &
        2&647 & 2&863 & 3&698 & 3&582  \\
\hline
 5&50 & \multicolumn{2}{|c|}{} & 1&523 & 1&108 & 0&053 & 0&100 &
        2&574 & 2&837 & 4&181 & 4&077  \\
\hline
 7&00 & \multicolumn{2}{|c|}{} & 2&197 & 1&608 & 0&053 & 0&100 &
       2&536 & 2&825 & 4&555 & 4&458  \\
\hline
 8&50 & \multicolumn{2}{|c|}{} & 2&872 & 2&108 & 0&052 & 0&100 &
       2&513 & 2&819 & 4&849 & 4&754  \\
\hline
 10&00 & \multicolumn{2}{|c|}{} & 3&548 & 2&608 & 0&050 & 0&100 &
       2&498 & 2&815 & 5&084 & 4&993  \\
\hline
\end{tabular}
\end{table}
Using a least squares fit, the following empirical formulae, related to the proposed method, are obtained from the results shown in Table 1, namely 
\begin{itemize}
\item optimum tuning points lying on $\Gamma_{y}$ for $0.1<=t_{p}<=0.7$
\begin{displaymath}
h=0.4541 - 0.1035 \, t_{p}+1.0794 \, t_{p}^{2}; \, h_{i}=0.8271 - 0.4805 \, t_{p}+0.5613 \, t_{p}^{2}
\end{displaymath}
\item optimum tuning points lying on $\Gamma_{v}$ for $0.85<=t_{p}<=10$
\begin{displaymath}
h=0.5884 + 0.5826 \, t_{p}+0.0033 \, t_{p}^{2}; \, h_{i}=0.7874 - 0.0434 \, t_{p}+0.0028 \, t_{p}^{2}
\end{displaymath}
\end{itemize}
The values of the parameters and performance indexes of the proposed method, as per these empirical formulae but considering all the decimal digits of the coefficients, are listed in Table 3. 

\begin{table}[htbp]
\caption{Empirical parameters and indexes of the proposed method}
\label{Tab3}
\centering
\begin{tabular}{|r @{.} l|r @{.} l|r @{.} l|r @{.} l|r @{.} l|r @{.} l|}
\hline
\multicolumn{2}{|c|}{$t_{p}$} &
\multicolumn{2}{|c|}{$h$} &
\multicolumn{2}{|c|}{$h_{i}$} &
\multicolumn{2}{|c|}{$PO_{y}$} &
\multicolumn{2}{|c|}{$PO_{v}$} &
\multicolumn{2}{|c|}{$ISE$} \\
\hline
 0&10 & 0&4546 & 0&7846 & 0&0096 & 0&0100 & 1&5259  \\
\hline
 0&25 & 0&4957 & 0&7420 & 0&0179 & 0&0361 & 1&6700  \\
\hline 
 0&40 & 0&5854 & 0&7247 & 0&0103 & 0&0512 & 1&7864  \\
\hline
 0&55 & 0&7237 & 0&7326 & 0&0096 & 0&0701 & 1&8762  \\
\hline
 0&70 & 0&9106 & 0&7657 & 0&0104 & 0&1052 & 1&9410  \\
\hline
\hline
 0&85 & 1&0861 & 0&7525 & \multicolumn{2}{|c|}{}& 0&0850 & 2&0509 \\
\hline
 1&00 & 1&1744 & 0&7468 & \multicolumn{2}{|c|}{}& 0&0974 & 2&1315 \\
\hline
 2&50 & 2&0658 & 0&6965 & \multicolumn{2}{|c|}{}& 0&1294 & 2&8900 \\
\hline
 4&00 & 2&9722 & 0&6589 & \multicolumn{2}{|c|}{}& 0&1115 & 3&5606 \\
\hline
 5&50 & 3&8935 & 0&6340 & \multicolumn{2}{|c|}{}& 0&0913 & 4&0928 \\
\hline
 7&00 & 4&8298 & 0&6218 & \multicolumn{2}{|c|}{}& 0&0821 & 4&4885 \\
\hline
 8&50 & 5&7810 & 0&6224 & \multicolumn{2}{|c|}{}& 0&0895 & 4&7722 \\
\hline
 10&00& 6&7473 & 0&6357 & \multicolumn{2}{|c|}{}& 0&1158 & 4&9695 \\
\hline
\end{tabular}
\end{table}

\section{Conclusion}
In this paper a new deadbeat response tuning method is applied to the PI controllers in first-order time-delay systems and compared with some existing ones. The involved performance indexes are strictly calculated by means of a new analytical solution of the differential difference equations relative to a step setpoint response. Moreover tuning charts, having as coordinates the two PI parameters  and provided with the stability zone borderline and with the performance indexes curves, have been introduced. These charts, previously prepared, may be considered a general and an effective tool, which does not require any advanced or theoretical knowledge, for the start-up and the maintenance of any controller.

\appendix{}

\section{Step setpoint response in systems without time delay}
During a step setpoint change from $r=0$ to $r=1$ applied to the closed-loop system in steady condition, the differential equation, relative to the controlled variable $y$ and deduced from (\ref{eq:1.4}), and its underdamped solution are given by
\begin{equation}\label{eq:A.1}
t_{i} \frac{d^{2}y(t)} {dt^{2}}\ + t_{i}(1+h) \frac{dy(t)} {dt}\ +h \, y(t)  = h  
\end{equation}
\begin{equation}\label{eq:A.2}
y(t)=1+e^{-a \,t}(-cos(b \, t)- \frac {a}{b}\ sin(b \,t))
\end{equation}
where
\begin{displaymath}
a=0.5(1+h) 
\end{displaymath}
\begin{equation}\label{eq:A.4}
b= \frac{0.5}{t_{i}}\ \left(4 \, h \, t_{i}-t_{i}^{2}(1+h)^{2} \right)^{0.5}
\end{equation}
For the controller output $u(t)$ multiplied by $K$, named $v(t)$, it follows from (\ref{eq:1.1}) and (\ref{eq:A.2})
\begin{equation}\label{eq:A.5}
\begin{split}
v(t)&=y(t)+ \frac {dy(t)}{dt}\ \\ &=1+e^{-a \,t}(-cos(b \, t)+ \frac {-a+a^{2} +b^{2}}{b}\ sin(b \, t))
\end{split}
\end{equation}
The derivatives of $y(t)$ and $v(t)$ with respect to $t$ are given respectively by
\begin{equation}\label{eq:A.6}
\frac{dy(t)}{dt}\ =e^{-a \,t} \frac {a^{2}+b^{2}}{b}\ sin(b \,t)
\end{equation}
\begin{equation}\label{eq:A.7}
\frac{dv(t)}{dt}\ = e^{-a \,t}(a^{2}+b^{2})(cos(b \, t)- \frac {a-1}{b}\ sin(b \, t))
\end{equation}
Introducing the solution $t$ of (\ref{eq:A.6}), equated to zero, in (\ref{eq:A.2}), we obtain for the overshoot $PO_{y}$ 
\begin{equation}\label{eq:A.8}
PO_{y} = e^{- \pi a/b}
\end{equation}
Analogously introducing the solution $t$ of (\ref{eq:A.7}), equated to zero, in (\ref{eq:A.5}), we obtain for the overshoot $PO_{v}$ 
\begin{equation}\label{eq:A.9}
PO_{v} = e^{ -(a/b)\arctan(b/(a-1))}(1-2 \,a+a^{2}+b^{2})^{0.5}
\end{equation}
From (\ref{eq:A.2}) it follows
\begin{equation}\label{eq:A.10}
ISE = \int_{0}^{\infty}(y(t)-1)^{2} \,dt = \frac{0.25}{a}\ \frac{5 \, a^{2}+b^{2}}{a^{2}+b^{2}}\
\end{equation}

\section{Step setpoint response in systems with time delay}
During a step setpoint change from $r=1$ to $r=0$ applied to the closed-loop system in steady condition, the differential equation, relative to the controlled variable $y$ and deduced from (\ref{eq:1.5}), is given by
\begin{equation}\label{eq:B.1}
\begin{split}
 \frac{dy(t)}{dt}\ &+ t_{p}  \frac {d^{2}y(t)}{dt^{2}}\ \\ &= (h_{i}+  h  \frac{d}{dt}\  ) (r(t-1)-y(t-1)) \quad  for \, t>=0
\end{split}
\end{equation}
subject to an initial condition of the form
\begin{equation}\label{eq:B.2}
 y(t)=y_{0}(t)=1 \quad for \quad -1<=t<=0
\end{equation}
The analytical solution of (\ref{eq:B.1}), obtained with the method of steps in  \cite{bib1},  consists of a set of functions $y_{n}(t_{n})$, each valid for $n-1<t<n$, whose independent variable $t_{n}$ has the temporal origin at $t=n-1$, given by   
\begin{equation}\label{eq:B.3}
\begin{split}
y_{1}(t_{1}) &= 1\\
y_{n}(t_{n}) &= \sum_{i=0}^{i=n-1}A_{n,i}t_{n}^{i}+ e^{-t_{n}/t_{p}} \sum_{j=0}^{j=n-2}B_{n,j}t_{n}^{j} \quad for \quad n>1
\end{split}
\end{equation}
In Example no.2 of  \cite{bib1} there are also the recursive expressions suitable for the evaluation of the coefficients $A_{n,i}$ and $B_{n,j}$.

For the controller output $u_{n}(t_{n})$ multiplied by $K$, named $v_{n}(t_{n})$, it follows from (\ref{eq:1.1})
\begin{equation}\label{eq:B.4}
v_{n}(t_{n})=y_{n+1}(t_{n+1})+ t_{p} \frac{dy_{n+1}(t_{n+1})}{dt_{n+1}}\
\end{equation}
Finally the integral of the squared error $ISE$ is calculated with
\begin{equation}\label{eq:B.5}
\begin{split}
ISE= &+ \frac{1}{2}\ 0.01 \,(y_{1}(0))^{2} - \frac{1}{2}\ 0.01 \, (y_{7}(1))^{2}\\ &+ 0.01 \sum_{n=1}^{n=7} \sum_{\tau=1}^{\tau=100} (y_{n} (0.01 \, \tau))^{2}\\ 
\end{split}
\end{equation}

\end{document}